\documentclass[final]{siamltex}
\usepackage{geometry}  
\usepackage{amsmath,amssymb,amsfonts,graphicx}
\usepackage{epsfig}
\usepackage{epstopdf}
\usepackage{subfigure, float}
\usepackage{algorithmic}
\usepackage{algorithm}
\usepackage{comment}
\usepackage[usenames,dvipsnames]{xcolor}
\newcommand{\citep}[1]{\cite{#1}}

\newcommand{\hdiv}{H(\mathrm{div})}
\newcommand{\hcurl}{H(\mathrm{curl})}

\newcommand{\dx}{d\mathbf{x}}

\title{Efficient discontinuous Galerkin finite element methods via Bernstein polynomials}
\author{Robert C. Kirby \footnotemark[2]}

\begin{document}
\maketitle

\pagestyle{myheadings}
\thispagestyle{plain}
\markboth{R.C. KIRBY}{Efficient DGFEM via Bernstein polynomials} 

\renewcommand{\thefootnote}{\fnsymbol{footnote}}

\footnotetext[2]{
Department of Mathematics, Baylor University; One Bear Place \#97328; Waco, TX 76798-7328.  This work is supported by NSF grant CCF-1325480.}

\renewcommand{\thefootnote}{\arabic{footnote}}

\begin{abstract}
We consider the discontinuous Galerkin method for hyperbolic
conservation laws, with some particular attention to the linear
acoustic equation, using Bernstein polynomials as local bases.
Adapting existing techniques leads to optimal-complexity
computation of the element and boundary flux terms.  The element mass
matrix, however, requires special care.  In particular, we give an explicit
formula for its eigenvalues and exact characterization of the
eigenspaces in terms of the Bernstein representation of orthogonal
polynomials.  We also show a fast algorithm for solving linear systems
involving the element mass matrix to preserve the overall complexity
of the DG method.  Finally, we present numerical results investigating
the accuracy of the mass inversion algorithms and the scaling of total
run-time for the function evaluation needed in DG time-stepping.
\end{abstract}



\section{Introduction}
Bernstein polynomials, which are ``geometrically decomposed'' in the
sense of~\citep{arnold2009geometric} and rotationally symmetric,
provide a flexible and general-purpose set of simplicial finite
element shape functions. 
Morever, recent research has demonstrated distinct algorithmic
advantages over other simplicial shape functions, as many essential
elementwise finite element computations can be performed on
with optimal complexity using Bernstein polynomials
In~\citep{kirby2011fast}, we showed
how, with constant coefficients, elementwise mass and stiffness
matrices could each be applied to vectors in $\mathcal{O}(n^{d+1})$ operations,
where $n$ is the degree of the local basis and $d$ is the spatial
dimension.  Similar blockwise linear algebraic structure enabled
quadrature-based algorithms in~\citep{kirby2012fast}.  Around the same time,
Ainsworth \emph{et al}~\citep{ainsworth2011bernstein} showed that the Duffy
transform~\citep{duffy1982quadrature} reveals a tensorial structure in
the Bernstein 
basis itself, leading to sum-factored algorithms for polynomial
evaluation and moment computation.  Moreover, they provide an
algorithm that assembles element matrices with $\mathcal{O}(1)$ work
per entry that utilizes their fast moment algorithm together with a
very special property of the Bernstein polynomials.  Work
in~\citep{kirbyssec,brown_sc_2012_20} extends these techniques to
$\hdiv$ and $\hcurl$. 

In this paper, we consider Bernstein polynomial techniques in a
different context -- discontinuous Galerkin methods for hyperbolic
conservation laws  
\begin{equation}
  \label{eq:conslaw}
q_t + \nabla \cdot F(q) = 0,
\end{equation}
posed on a domain 
$\Omega \times [0,T) \subset \mathbb{R}^d \times\mathbb{R}$, together
with suitable initial and boundary conditions.  As a
particular example, we consider the linear acoustic model
\begin{equation}
\begin{split}
p_{t} + \nabla \cdot u & = 0,\\
u_{t} + \nabla p & = 0,
\end{split}
\end{equation}
Here, $q = [u,p]^T$ where the pressure variable $p$ is a scalar-valued function on
$\Omega \times [0,T]$ and the velocity $u$ maps the same space-time
domain into $\mathbb{R}^d$.  

Discontinuous Galerkin (DG) methods for such problems place finite
volume methods in a variational 
framework and extend them to higher orders of polynomial
approximation~\cite{cockburn1991runge}, but fully realizing the potential
efficiencies of high-order methods requires careful consideration of
algorithmic issues.
Simplicial orthogonal polynomials~\cite{dubiner1991spectral,KarShe05}
provide one existing mechanism for achieving low operation counts.  Their
orthogonality gives diagonal local mass 
matrices.  Optimality then requires special quadrature that reflects the
tensorial nature of the basis under the Duffy transform or
collapsed-coordinate mapping from the $d$-simplex to the $d$-cube and
also includes appropriate points to incorporate contributions
from both volume and boundary flux terms.
Hesthaven and Warburton~\citep{hesthaven2002nodal,hesthaven2007nodal}
propose an alternate approach, 
using dense linear algebra in conjunction with Lagrange polynomials.
While of greater algorithmic complexity, highly-tuned
matrix multiplication can make this approach
competitive or even superior at practical polynomial orders.
Additional extensions of this idea include the so-called ``strong DG''
forms and also a pre-elimination of the elementwise mass matrix giving
rise to a simple ODE system.  With care, this approach can give very
high performance on both CPU and GPU systems~\citep{klockner2009nodal}.

In this paper, we will show how each term in the DG
formulation with Bernstein polynomials as the local basis can be
handled with optimal complexity   For the element 
and boundary flux terms, this requires only an adaptation of existing
techniques, but inverting the element mass matrix turns out to be a
challenge lest it dominate the complexity of the entire process.
We rely on the recursive block structure
described in~\citep{kirby2011fast} to give an $\mathcal{O}(n^{d+1})$
algorithm for solving linear systems with the constant-coefficient
mass matrix.  
We may view our approach as sharing certain
important features of both collapsed-coordinate
and Lagrange bases.  Like collapsed-coordinate methods, we seek
to use specialized  structure to optimize algorithmic complexity.
Like Lagrange polynomials, we seek to do this using a relatively
discretization-neutral basis.

\section{Discontinuous Galerkin methods}
We let $\mathcal{T}_h$ be a triangulation of $\Omega$ in the sense of
~\citep{BreSco} into affine simplices.  For curved-sided elements, we could
adapt the techniques of~\cite{warburton2013low} to incorporate the Jacobian into our local basis functions to recover the reference mass matrix on each
cell at the expense of having variable coefficients in other
operators, but this does not affect the overall order of complexity.
We let $\mathcal{E}_h$ denote the set of all edges in the triangulation.

For $T \in \mathcal{T}_h$, let $P_n(T)$ be the space of polynomials of
degree no greater than $n$ on $T$.  This is a vector space of
dimension  $P^d_n \equiv \binom{n+d}{n}$.  We define the global finite
element space
\begin{equation}
V_h = \left\{ f : \Omega \rightarrow \mathbb{R} : f|_T \in P_n(T), \ T
\in \mathcal{T}_h \right\},
\end{equation}
with no continuity
enforced between cells.  
Let $\left( \cdot , \cdot \right)_T$ denote the standard $L^2$ inner
product over $T \in \mathcal{T}_h$ and 
$\langle \cdot , \cdot \rangle_\gamma$ the $L_2$ 
inner product over an edge $\gamma \in \mathcal{E}_h$.

After multiplying~\eqref{eq:conslaw} by a test function and integating by parts elementwise, a DG method seeks $u_h$ in $V_h$ such that
\begin{equation}
\label{eq:dg}
\sum_{T \in \mathcal{T}_h} 
\left[
\left( u_{h,t} , v_h \right)_T
- \left( F(u_h) , \nabla v_h \right)_T \right]
+\sum_{\gamma \in \mathcal{E}_h} \langle \hat{F} \cdot n , v_h \rangle
\end{equation}
for all $v_h \in V_h$.

Fully specifying the DG method requires defining a numerical flux
function $\hat{F}$ on each $\gamma$.  On internal edges, it takes
values from either side of the 
edge and produces a suitable approximation to the flux $F$.  Many
Riemann solvers from the finite volume literature have
been adapted for DG
methods~\cite{cockburn1991runge,dumbser2008unified,toro1999riemann}.
The particular choice of numerical flux does not matter for our
purposes.  On external edges,  
we choose $\hat{F}$ to appropriately enforce boundary conditions.

This discretization gives rise to a system of ordinary differential
equations 
\begin{equation}
\mathrm{M} \mathrm{u}_t + \mathrm{F}(\mathrm{u}) = 0,
\end{equation}
where $\mathrm{M}$ is the block-diagonal mass matrix and
$\mathrm{F}(\mathrm{u})$ includes the cell and boundary flux terms.
Because of the hyperbolic nature of the system, explicit methods are
frequently preferred.  A forward Euler method, for example, gives
\begin{equation}
\mathrm{u}^{n+1} = \mathrm{u}^n - \Delta t \mathrm{M}^{-1}
\mathrm{F}(\mathrm{u^n}) \equiv \mathrm{u}^n - \Delta t L(\mathrm{u}^n),
\end{equation}
which requires the application of $\mathrm{M}$ at each time step.
The SSP methods~\cite{gottlieb2001strong,shu1988total} give stable
higher-order in time methods.  For example, the well-known third order
scheme is \begin{equation}
\begin{split}
\mathrm{u}^{n,1} & = \mathrm{u}^n + \Delta t L(\mathrm{u}^n), \\
\mathrm{u}^{n,2} & = \frac{3}{4} \mathrm{u}^n + \frac{1}{4}
\mathrm{u}^{n,1} + \frac{1}{4} \Delta t L(\mathrm{u}^{n,1}),\\
\mathrm{u}^{n+1} & = \frac{1}{3} \mathrm{u}^n + \frac{2}{3}
\mathrm{u}^{n,2} + \frac{2}{3} \Delta t L( \mathrm{u}^{n,2}).
\end{split}
\end{equation}

Since the Bernstein polynomials give a dense element mass
matrix, applying $\mathrm{M}^{-1}$ efficiently will
require some care.  It turns out that $\mathrm{M}$ possesses many
fascinating properties that we shall survey in Section~\ref{sec:mass}.
Among these, we will give an $\mathcal{O}(n^{d+1})$ algorithm for
applying the elementwise inverse.

DG methods yield reasonable solutions to acoustic or Maxwell's
equations without slope limiters, although most nonlinear problems will
require them to suppress oscillations.  Even linear transport can
require limiting when a discrete maximum principle is
required. Limiting high-order polynomials on simplicial domains
remains quite a challenge.  It may be possible to utilize properties
of the Bernstein polynomials to design new limiters or conveniently
implement existing ones. 
For example, the convex hull property (i.e. that polynomials in the
Bernstein basis lie in the convex hull of their control points) gives
sufficient conditions for enforcing extremal bounds.
We will not offer further contributions in this direction, but refer
the reader to other works on higher order limiting such
as~\citep{hoteit2004new,zhu2008runge,zhu2013runge}. 



\section{Bernstein-basis finite element algorithms}
\label{BBFEA}

\subsection{Notation for Bernstein polynomials}
We formulate Bernstein polynomials on the $d$-simplex 
using barycentric coordinates and multiindex notation.  For a
nondegenerate simplex $T \subset \mathbb{R}^{d}$ with vertices
$\{x_i\}_{i=0}^{d}$, let $\{b_i\}_{i=0}^{d}$ denote the
barycentric coordinates.  Each $b_i$ affinely maps
$\mathbb{R}^d$ into $\mathbb{R}$ with $b_i(x_j) = \delta_{ij}$ for 
$0 \leq i,j \leq d$.  It follows that $b_i(x) \geq 0$ for all $x \in T$.

We will use common multiindex notation, denoting multiindices
with Greek letters, although we will begin the indexing with 0 rather than 1.
So,
$\alpha=(\alpha_0,\alpha_1,\dots,\alpha_{d})$ is a tuple of $d+1$
nonnegative integers.  We define the \emph{order} of a multiindex $\alpha$ 
by $\left| \alpha \right| \equiv \sum_{i=0}^{d} \alpha_i$.  We say
that $\alpha \geq \beta$ provided that the inequality 
$\alpha_i \geq \beta_i$ holds componentwise for $0 \leq i \leq d$.  
Factorials and binomial coefficients over multiindices have implied multiplication.  That is,
\[ 
\alpha! \equiv \prod_{i=0}^{d} \alpha_i!
\]
and, provided that $\alpha \geq \beta$,
\[
\binom{\alpha}{\beta} = \prod_{i=0}^d \binom{\alpha_i}{\beta_i}.
\]
Without ambiguity of notation, we also define a binomial coefficient
with a whole number for the upper argument and and multiindex as lower
by
\[
\binom{n}{\alpha} = \frac{n!}{\alpha!} = \frac{n!}{\prod_{i=0}^n \alpha_i!}.
\]
We also define $e_i$ to be the multiindex consisting of zeros in all but the $i^{\mathrm{th}}$ entry, where it is one.  

Let $\mathbf{b} \equiv \left( b_0 , b_2 , \dots , b_d \right) $ be
a tuple of barycentric coordinates on a simplex.  For multiindex
$\alpha$, we define a \emph{barycentric monomial} by
\[
\mathbf{b}^\alpha = \prod_{i=0}^{d} b_i^\alpha.
\]
We obtain the Bernstein polynomials by scaling these by certain
binomial coefficients  
\begin{equation}
B^n_\alpha =  \frac{n!}{\alpha!} \mathbf{b}^\alpha.
\end{equation}
For all spatial dimensions and degrees $n$, the Bernstein
polynomials of degree $n$
\[
\left\{
B^n_{\alpha}
\right\}_{\left| \alpha \right| = n},
\]
form a nonnegative partition of unity and a basis for the vector space of
polynomials of degree $n$.  They are suitable for
assembly in a $C^0$ fashion or even into smoother splines~\cite{LaiSch07}.
While DG methods do not require assembly, the geometric decomposition
does make handling the boundary terms straightforward.

Crucial to fast algorithms using the Bernstein basis, as
originally applied to $C^0$
elements~\citep{ainsworth2011bernstein,kirby2011fast}, is the 
sparsity of differentiation.  That is, it takes no more than $d+1$
Bernstein polynomials of degree $n-1$ to represent the derivative of a
Bernstein polynomial of degree $n$.

For some coordinate direction $s$, we use the general product rule to write 
\[
\frac{\partial B^n_\alpha}{\partial s}
= \frac{\partial}{\partial s} \left( \frac{n!}{\alpha!} \mathbf{b}^\alpha \right)
= \frac{n!}{\alpha!} \sum_{i=0}^{d} \left( \alpha_i \frac{\partial b_i}{\partial s} b_i^{\alpha_i-1} \Pi_{j=0}^{d} b_i^{\alpha_i} \right),
\]
with the understanding that a term in the sum vanishes if $\alpha_i = 0$.  This can readily be rewritten as
\begin{equation}
	\frac{\partial B^n_\alpha}{\partial s} 
	= n \sum_{i=0}^{d}  B_{\alpha-e_i}^{n-1} \frac{\partial b_i}{\partial s},
\end{equation}
again with the terms vanishing if any $\alpha_i = 0$,
so that the derivative of each Bernstein polynomial is a short linear combination of lower-degree Bernstein polynomials.

Iterating over spatial directions, the gradient of each Bernstein polynomial can be written as
\begin{equation}
	\label{eq:berngrad}
	\nabla  B^n_\alpha 
	= n \sum_{i=0}^{d} B_{\alpha-e_i}^{n-1} \nabla b_i.
\end{equation}
Note that each $\nabla b_i$ is a fixed vector in $\mathbb{R}^n$ for a
given simplex $T$.  In~\cite{kirbyssec}, we provide a data structure
called a \emph{pattern} for representing gradients as well as exterior
calculus basis functions.  For implementation details, we refer the
reader back to~\cite{kirbyssec}.

The \emph{degree elevation} operator will also play a crucial role in
our algorithms.  This operator expresses a B-form polynomial of degree $n-1$
as a degree $n$ polynomial in B-form.  For the orthogonal and
hierarchical bases in~\cite{KarShe05}, this operation would be trivial
-- appending the requisite number of zeros in a vector, while for
Lagrange bases it is typically quite dense.  
Whiel not trivial, degree elevation for Bernstein polynomials is still
efficient.  Take any Bernstein polynomial and multiply it by
$\sum_{i=0}^db_i = 1$ to find
\begin{equation}
\begin{split}
B^{n-1}_\alpha & = \left( \sum_{i=0}^{d} b_i \right) B^{n-1}_\alpha 
 = \sum_{i=0}^d b_i B^{n-1}_\alpha \\
& = \sum_{i=0}^d \frac{(n-1)!}{\alpha!} \mathbf{b}^{\alpha+e_i}
= \sum_{i=0}^d \frac{\alpha_i+1}{n} \frac{n!}{\left(\alpha+e_i\right)!}
\mathbf{b}^{\alpha+e_i} \\
& = \sum_{i=0}^d \frac{\alpha_i+1}{n} B^n_{\alpha+e_i}.
\end{split}
\end{equation}
We could encode this operation as a $P^d_n \times P^d_{n-1}$
matrix consisting of exactly $d+1$ nonzero entries, but it can also
be applied with a simple nested loop.  At any rate, we denote this
linear operator as $E^{d,n}$, where $n$ is the degree of the
resulting polynomial.  We also denote 
$E^{d,n_1,n_2}$ the operation that successively raises a polynomial
from degree $n_1$ into $n_2$.  This is just the product of $n_2 - n_1$
(sparse) operators:
\begin{equation}
E^{d,n_1,n_2} = E^{d,n_2} \dots E^{d,n_1+1}.
\end{equation}
We have that $E^{d,n} = E^{d,n-1,n}$ as a special case.

\subsection{Stroud conical rules and the Duffy transform}
The Duffy transform~\cite{duffy1982quadrature} tensorializes the
Bernstein polynomials, so sum factorization can be used for
evaluating and integrating these polynomials with Stroud conical
quadrature.  We used similar quadrature rules in our own work on
Bernstein-Vandermonde-Gauss matrices~\cite{kirby2012fast}, but the connection to
the Duffy transform and decomposition of Bernstein polynomials was
quite cleanly presented by Ainsworth \emph{et al} in~\cite{ainsworth2011bernstein}. 

The Duffy transform maps any point $\mathbf{t} = (t_1,t_2,\dots,t_n)$
in the $d$-cube $[0,1]^n$ into the barycentric coordinates for a
$d$-simplex by first defining 
\begin{equation}
	\lambda_0 = t_1
\end{equation}
and then inductively by
\begin{equation}
	\lambda_{i} = t_{i+1} \left( 1 - \sum_{j=0}^{i-1} \lambda_j \right)
\end{equation}
for $1 \leq i \leq d-1$, and then finally
\begin{equation}
	\lambda_{d} = 1 - \sum_{j=0}^{d-1} \lambda_j.
\end{equation}

If a simplex $T$ has vertices $\{ \mathbf{x}_i \}_{i=0}^{d}$, then the mapping
\begin{equation}
	\mathbf{x}(\mathbf{t}) = \sum_{i=0}^{d} \mathbf{x}_i \lambda_i(\mathbf{t})
\end{equation}
maps the unit $d$-cube onto $T$.

This mapping can be used to write integrals over $T$ as iterated
weighted integrals over $[0,1]^d$ 
\begin{equation}
	\int_T f(\mathbf{x}) d\mathbf{x} = \frac{|T|}{d!} \int_0^1
        dt_1 (1-t_1)^{d-1} \int_0^1 dt_2 (1-t_2)^{d-2} \dots \int_0^1
        dt_t f(\mathbf{x}(t)). 
\end{equation} 
The \emph{Stroud conical rule}~\cite{stroud} is based on this observation
and consists of tensor products of certain Gauss-Jacobi quadrature
weights in each $t_i$ variable, where the weights are chosen to absorb
the factors of $(1-t_i)^{n-i}$.  These rules play an important role in
the collapsed-coordinate framework of~\cite{KarShe05} among many other
places.

As proven in~\cite{ainsworth2011bernstein}, pulling the Bernstein
basis back to $[0,1]^d$ under the Duffy transform reveals a
tensor-like structure.
It is shown that with $B^n_i(t) = \binom{n}{i} t^i (1-t)^{n-i}$ the
one-dimensional Bernstein polynomial, that  
\begin{equation}
	B^n_\alpha(\mathbf{x}(\mathbf{t})) = B_{\alpha_0}^n(t_1) B_{\alpha_1}^{n-\alpha_0}(t_2) \cdots 
	B_{\alpha_{d-1}}^{n-\sum_{i=0}^{d-2}\alpha_i}(t_d).
\end{equation}
This is a ``ragged'' rather than true tensor product, much as the
collapsed coordinate simplicial bases~\cite{KarShe05}, but entirely
sufficient to enable sum-factored algorithms.

\subsection{Basic algorithms}
The Stroud conical rule and tensorialization of Bernstein polynomials
under the Duffy transformation lead to highly efficient algorithms for
evaluating B-form polynomials and approximating moments of functions
against sets of Bernstein polynomials. 

Three algorithms based on this decomposition turns out to be fundamental for
optimal assembly and application of Bernstein-basis bilinear forms.
First, any polynomial 
$u(\mathbf{x})=\sum_{|\alpha|=n} \mathrm{u}_\alpha
B^n_\alpha(\mathbf{x})$ 
may be evaluated at the Stroud conical points in
$\mathcal{O}(n^{d+1})$ operations.  In~\cite{kirby2012fast}, this
result is presented as exploiting certain block structure in the
matrix tabulating the Bernstein polynomials at quadrature points.
In~\cite{ainsworth2011bernstein}, it is done by explicitly factoring
the sums. 

Second, given some function $f(\mathbf{x})$ tabulated at the Stroud points, it
is possible to approximate the set of Bernstein moments 
\[
\mu^n_\alpha(f) = \int_T f(\mathbf{x}) B^n_\alpha \dx
\]
for all $|\alpha|=n$
via Stroud quadrature in $\mathcal{O}(n^{d+1})$ operations.  
In the the case where $f$ is constant on
$T$, we may also use the algorithm for applying a mass matrix
in~\cite{kirby2011fast} to bypass numerical integration. 

Finally, it is shown in~\cite{ainsworth2011bernstein} that the moment
calculation can be adapted to the evaluation of element mass and hence
stiffness and convection matrices utilizing another remarkable
property of the Bernstein polynomials.  Namely, the product of two
Bernstein polynomials of any degrees is, up to scaling, a Bernstein
polynomial of higher degree:
\begin{equation}
	B_\alpha^{n_1} B_\beta^{n_2} =
        \frac{\binom{\alpha+\beta}{\alpha}}{\binom{n_1+n_2}{n_1}}
        B^{n_1+n_2}_{\alpha+\beta}, 
\end{equation}

Also, the first two algorithms described above for evaluation and
moment calculations demonstrate that $M$ may be applied to a vector
without explicitly forming its entries in only $\mathcal{O}(n^{d+1})$
entries.  In~\cite{kirbyssec}, we show how to adapt these algorithms to short linear
combinations of Bernstein polynomials so that stiffness and convection
matrices require the same  order of complexity as the mass.

\subsection{Application to DG methods}
As part of each explicit time stepping stage, we must evaluate
$\mathrm{M}^{-1} \mathrm{F(u)}$.  Evaluating $\mathrm{F(u)}$ requires
handling the two flux terms in~\eqref{eq:dg}.  To handle
\[
\left( F(u_h) , \nabla v_h \right)_T, 
\]
we simply evaluate $u_h$ at the Stroud points on $T$, which requires
$\mathcal{O}(n^{d+1})$ operations.  Then, evaluating $F$ at each of
these points is purely pointwise and so requires but
$\mathcal{O}(n^{d})$.  Finally, the moments against gradients of
Bernstein polynomials also requires $\mathcal{O}(n^{d+1})$
operations.  This term, then, is readily handled by existing Bernstein
polynomial techniques.

Second, we must address, on each interface $\gamma \in \mathcal{E}$,
\[
\langle \hat{F} \cdot n , v_h \rangle_\gamma.
\]
The numerical flux $\hat{F} \cdot n$ requires the values of $u_h$ on
each side of the interface and is evaluated pointwise at each facet
quadrature point.  Because of the Bernstein polynomials'
geometric decomposition, only $P^{d-1}_n$ basis functions are nonzero
on that facet, and their traces are in fact exactly the Bernstein
polynomials on the facet.  So we have to evaluate two polynomials (the
traces from each side) of
degree $n$ in $d-1$ variables at the facet Stroud points.  This
requires $\mathcal{O}(n^d)$ operations.  The numerical flux is
computed pointwise at the $\mathcal{O}(n^{d-1})$ points, and then the
moment integration is performed on facets for an overall cost of
$\mathcal{O}(n^d)$ for the facet flux term.
In fact, the geometric decomposition makes this term much easier to
handle optimally with Bernstein polynomials than collapsed-coordinate
bases, although though specially adapted Radau-like quadrature rules,
the boundary sums may be lifted into the volumetric
integration~\citep{tim}.

The mass matrix, on the other hand, presents a much deeper challenge
for Bernstein polynomials than for collapsed-coordinate ones.  Since
it is dense with 
$\mathcal{O}(n^d)$ rows and columns, a standard matrix
Cholesky decomposition requires $\mathcal{O}(n^{3d})$
operations as a startup cost, followed by a pair of triangular solves 
on each solve at $\mathcal{O}(n^{2d})$ each. 
For $d>1$, this complexity clearly dominates the steps
above, although an optimized Cholesky routine might very well win at
practical orders.   In the next section, we turn to a careful study of
the mass matrix, deriving an algorithm of optimal complexity. 

\section{The Bernstein mass matrix}
\label{sec:mass}
We begin by defining the rectangular
Bernstein mass matrix on a $d$-simplex $T$ by
\begin{equation}
M^{T,m,n}_{\alpha\beta} = \int_T B^m_{\alpha} B^n_{\beta} dx,
\end{equation}
where $m,n \geq 0$.  

By a change of variables, we can write
\begin{equation}
M^{T,m,n} = M^{d,m,n} |T| d!,
\end{equation}
where $M^{d,m,n}$ is the mass matrix on the unit right simplex $S_d$ in
$d$-space and $|T|$ is the $d$-dimensional measure of $T$. 
When $m=n$, we suppress the third superscript and write $M^{T,m}$ or $M^{d,m}$.
We include the more general case of a
rectangular matrix because such will appear later in our discussion of the
block structure.

This mass matrix has many beautiful properties.  Besides the block-recursive structure developed in~\citep{kirby2011fast},
it is related to the Bernstein-Durrmeyer
operator~\cite{derriennic1985multivariate,farouki2003construction} of approximation theory.  Via this connection, we 
provide an exact characterization of its eigenvalues and
associated eigenspaces in the square case $m=n$.  Finally, and most
pertinent to the case of discontinuous Galerkin methods, we describe
algorithms for solving linear systems involving the mass matrix.

Before proceeding, we recall from~\cite{kirby2011fast} that,  
formulae for integrals of products of powers of barycentric
coordinates, the mass matrix formula is exactly
\begin{equation}
\label{eq:Mdef}
M^{d,m,n}_{\alpha,\beta} = \frac{m!n!\left( \alpha +
  \beta\right)!}{\left(m+n+d\right)!\alpha!\beta!}
\end{equation}

\subsection{Spectrum}
The Bernstein-Durrmeyer
operator~\cite{derriennic1985multivariate} is defined on $L^2$ by
\begin{equation}
D_n(f) = \frac{\left(n+d\right)!}{n!} \sum_{|\alpha|=n} \left( f , B^n_\alpha \right).
\end{equation}
This has a structure similar to a discrete Fourier series, although
the Bernstein polynomials are orthogonal.  The original Bernstein
operator~\citep{LaiSch07} has the form of a Lagrange interpolant,
although the basis is not interpolatory.

For $i\geq 1$, we let $Q_i$ denote the space of $d$-variate
polynomials of degree $i$ that are $L^2$ orthogonal to all polynomials of degree $i-1$ on the simplex.
The following result is given in~\cite{derriennic1985multivariate},
and also referenced in~\cite{farouki2003construction} to generate the
B-form of simplicial orthogonal polynomials.
\begin{theorem}[Derriennic]
For each $0 \leq i \leq n$, each 
\[
\lambda_{i,n} = \frac{\left(n+d\right)! n!}{\left( n + i + d \right)!
  \left( n - i \right)!}
\]
is an eigenvalue of $D_n$ corresponding to the eigenspace $Q_i$.
\end{theorem}

This gives a sequence of eigenvalues 
$\lambda_{0,n} > \lambda_{1,n} > \dots > \lambda_{n,n} > 0$, each
corresponding to polynomial eigenfunctions of increasing degree.

Up to scaling, the Bernstein-Durrmeyer operator restricted to
polynomials $P_n$ exactly corresponds to the action of the mass
matrix.  To see this, suppose that $P_n \ni p = \sum_{|\alpha|=n}
\mathrm{p}_\alpha B^n_\alpha$.  Then
\begin{equation}
\begin{split}
\frac{n!}{\left(n+d\right)!} D_n(p) & = \sum_{|\alpha|=n} \left( p ,
B^n_\alpha \right) B^n_\alpha \\
& = \sum_{|\alpha|=n} \left( \sum_{|\beta|=n} p_\beta B^n_\beta ,
B^n_\alpha \right) B^n_\alpha \\
& = \sum_{|\alpha|=n} \sum_{|\beta|=n} p_\beta \left( B^n_\beta ,
B^n_\alpha \right) B^n_\alpha \\
& = \sum_{|\alpha|=n} \left( \sum_{|\beta|=n} M^n_{\alpha,\beta}
p_\beta  \right) B^n_\alpha \\
\end{split}
\end{equation}
This shows that the coefficients of the B-form of 
$\frac{n!}{\left(n+d\right)!} D_n(p)$ are just the entries of the
Bernstein mass matrix times the coefficients of $p$.  Consequently,
\begin{theorem}
  \label{thm:eigs}
For each $0 \leq i \leq n$, each 
\[
\lambda_{i,n} = \frac{\left(n+d\right)! \left(n!\right)^2}{\left( n + i + d \right)!
  \left( n - i \right)!}
\]
is an eigenvalue of $M^n$ of multiplicity of $\binom{d+i-1}{d-1}$, and
the eigenspace is spanned by the B-form of any basis for $Q_i$.
\end{theorem}

This also implies that the Bernstein mass matrices are quite
ill-conditioned in the two norm, using the characterization in terms
of extremal eigenvalues for SPD matrices.
\begin{corollary}
  \label{cor:badguy}
The 2-norm condition number of $M^{d,n}$ is 
\begin{equation}
\frac{\lambda_{0,n}}{\lambda_{n,n}}
= \frac{(2n+d)!}{(n+d)!n!}
\end{equation}
\end{corollary}

However, the spread in eigenvalues does not tell the whole story.
We have exactly $n+1$ distinct eigenvalues,
independent of the spatial dimension.  This shows significant
clustering of eigenvalues when $d \geq 1$.
\begin{corollary}
  \label{cor:wonthappen}
In exact arithmetic, unpreconditioned conjugate gradient iteration
will solve a linear system of the form $M^{d,n} x = y$ in exactly
$n+1$ iterations, independent of $d$.
\end{corollary}

If the fast matrix-vector algorithms
in~\citep{ainsworth2011bernstein,kirby2011fast} 
are used to compute the matrix-vector product, this gives a total 
operation count of $\mathcal{O}(n^{d+2})$.  Interestingly, this ties
the per-element cost of Cholesky factorization when $d=2$, but without
the startup or storage cost.  It even beats a pre-factored matrix when
$d \geq 2$, but still loses asymptotically to the cost of
evaluating $\mathrm{F(u)}$.
However, in light of the large condition number given by
Corollary~\ref{cor:badguy}, it is doubtful whether this iteration
count can be realized in actual floating point arithmetic.

The high condition number also suggests an additional source of error
beyond discretization error.  Suppose that we commit an error of order
$\epsilon$ in solving $M x = y$, computing instead some $\hat{x}$ such
that $\left\| x - \hat{x} \right\| = \epsilon$ in the $\infty$ norm.
Let $u$ and $\hat{u}$ be the polynomial with B-form coefficients $x$
and $\hat{x}$, respectively.  Because a polynomial in B-form lies in
the convex hull of its control points~\citep{LaiSch07}, we also know
that $u$ and $\hat{u}$ differ by at most this same $\epsilon$ in the
max-norm.  Consequently, the roundoff error in mass inversion can
conceivably pollute the finite element approximation at high order,
although ten-digit accuracy, say, will still only give a maximum of $10^{-10}$
additional pointwise error in the finite element solution -- typically
well below discretization error.

\subsection{Block structure and a fast solution algorithm}
Here, we recall several facts proved in~\citep{kirby2011fast} related
to the block structure of $M^{d,m,n}$, which we will apply now for
solving square systems.

We consider partitioning the mass matrix formula~(\ref{eq:Mdef}) by
freezing the first entry in $\alpha$ and $\beta$.  Since there are
$m+1$ possible values for 
for $\alpha_0$ and $n+1$ for $\beta_0$, this partitions
$M^{d,m,n}$ into an $(m+1) \times (n+1)$ array, with blocks
of varying size.  In fact, each block
$M^{d,m,n}_{\alpha_0,\beta_0}$ is 
$P^{d-1}_{m-\alpha_0} \times P^{d-1}_{n-\beta_0}$.

These blocks are themselves, up to scaling,
Bernstein mass matrices of lower dimension.  In particular, we showed that
\begin{equation}
\label{eq:Mblock}
M^{d,m,n}_{\alpha_0,\beta_0} = 
\frac{\binom{m}{\alpha_0}
  \binom{n}{\beta_0}}{\binom{m+n+d+1}{\alpha_0+\beta_0} \left( m + n +
  d \right)} 
M^{d-1,m-\alpha_0,n-\beta_0}.
\end{equation}
We introduce the $(m+1)\times (n+1)$ array consisting of the scalars
multiplying the lower-dimensional mass matrices as
\begin{equation}
  \label{eq:nu}
\nu^{d,m,n}_{\alpha_0,\beta_0} = 
\frac{\binom{m}{\alpha_0}
  \binom{n}{\beta_0}}{\binom{m+n+d+1}{\alpha_0+\beta_0} \left( m + n +
  d \right)} 
\end{equation}
so that $M^{d,m,n}$ satisfies the block structure, with superscripts
on $\nu$ terms dropped for clarity
\begin{equation}
M^{d,m,n} = 
\begin{pmatrix}
\nu_{0,0} M^{d-1,m,n} & \nu_{0,1} M^{d-1,m,n-1} & \dots & \nu_{0,n} M^{d-1,m,0} \\
\nu_{1,0} M^{d-1,m-1,n} & \nu_{1,1} M^{d-1,m-1,n-1} & \dots &
\nu_{1,n} M^{d-1,m-1,0} \\
\vdots & \vdots & \ddots & \vdots \\
\nu_{n,0} M^{d-1,0,n} & \nu_{n,1} M^{d-1,0,n-1} & \dots & \nu_{n,n} M^{d-1,0,0} \\
\end{pmatrix}.
\end{equation}

We partition the right-hand side and solution vectors $y$ and $x$
conformally to $M$, so that the block $y_j$ is of dimension
$P^{d-1}_{n-j}$ and corresponds to a polynomial's $B$-form
coefficients with first indices equal to $j$.
We write the linear system in an augmented block matrix as
\begin{equation}
\label{eq:aug}
\left(
\begin{array}{cccc|c}
\nu_{0,0} M^{d-1,n,n} & \nu_{0,1} M^{d-1,n,n-1} & \dots & \nu_{0,n} M^{d-1,n,0} & y_0 \\
\nu_{1,0} M^{d-1,n-1,n} & \nu_{1,1} M^{d-1,n-1,n-1} & \dots & \nu_{1,n} M^{d-1,n-1,0} &
y_1 \\
\vdots & \vdots & \ddots & \vdots \\
\nu_{n,0} M^{d-1,n-1,n} & \nu_{n,1} M^{d-1,n-1,n-1} & \dots & \nu_{n,n} M^{d-1,0,0} & y_n \\
\end{array}
\right).
\end{equation}

From~\citep{kirby2011fast}, we also know that mass matrices of the
same dimension but differing degrees are related via degree elevation
operators by 
\begin{equation}
\label{eq:elTM}
M^{d,m-1,n} = \left( E^{d,m} \right)^t M^{d,m,n}.
\end{equation}
and
\begin{equation}
\label{eq:Mel}
M^{d,m,n-1} = M^{d,m,n} E^{d,n}.
\end{equation}
Iteratively, these results give
\begin{equation}
  \label{eq:elevateforward}
M^{d,m-i,n} = \left( E^{d,m-i,m} \right)^T  M^{d,m,n}.
\end{equation}
for $1 \leq i \leq m$ and
\begin{equation}
  \label{eq:elevatebackward}
M^{d,m,n-j} = M^{d,m,n} E^{d,n-j,n}
\end{equation}
for $1 \leq j \leq n$.
In~\citep{kirby2011fast}, we used these features to provide a fast
algorithm for matrix multiplication, but here we use them to
efficiently solve linear systems.

Carrying out blockwise Gaussian elimination in~\eqref{eq:aug}, we multiply the first row,
labeled with 0 rather than 1, by 
$\frac{\nu_{1,0}}{\nu_{0,0}} M^{d-1,n-1,n} \left( M^{d-1,n,n} \right)^{-1}$ 
and subtract from row 1 to introduce a zero
block below the diagonal.  However, this simplifies, as~\eqref{eq:elTM} tells us that
\begin{equation}
\label{eq:elimMtrick}
M^{d-1,n-1,n} 
\left(M^{d-1,n,n} \right)^{-1}
= \left( E^{d-1,n} \right)^t M^{d-1,n,n} \left(M^{d-1,n,n}
\right)^{-1}
= \left( E^{d-1,n} \right)^t.
\end{equation}
Because of this, along row 1 for $j \geq 1$, the elimination step gives entries of the form
\[
\nu_{1j} M^{d-1,n-1,n-j} - \frac{\nu_{10}\nu_{0j}}{\nu_{00}} \left( E^{d-1,n} \right)^t M^{d-1,n,n-j},
\]
but~\eqref{eq:elTM} renders this as simply
\begin{equation}
\nu_{1j} M^{d-1,n-1,n-j} - \frac{\nu_{10}\nu_{0j}}{\nu_{00}} M^{d-1,n-1,n-j}
= \left( \nu_{1j} - \frac{\nu_{10}\nu_{0j}}{\nu_{00}} \right) M^{d-1,n-1,n-j}.
\end{equation}
That is, the row obtained by block Gaussian elimination is the same as
one would obtain simply by performing a step of Gaussian elimination on the
matrix of coefficients $N^{d,n}$ containing the $\nu$ values above, as the
matrices those coefficients scale do not change under the row
operations.  Hence, performing elimination on the $(n+1)\times(n+1)$
matrix, independent of the dimension $d$, forms a critical step in the
elimination process.  After the block upper triangularization, we arrive at
a system of the form

\begin{equation}
\label{eq:augut}
\left(
\begin{array}{cccc|c}
\widetilde{\nu}_{0,0} M^{d-1,n,n} & \widetilde{\nu}_{0,1} M^{d-1,n,n-1} & \dots & \widetilde{\nu}_{0,n} M^{d-1,n,0} & \widetilde{y}_0 \\
0 & \widetilde{\nu}_{1,1} M^{d-1,n-1,n-1} & \dots & \widetilde{\nu}_{1,n} M^{d-1,n-1,0} &
\widetilde{y}_1 \\
\vdots & \vdots & \ddots & \vdots \\
0 & 0 & \dots & \widetilde{\nu}_{n,n} M^{d-1,0,0} & \widetilde{y}_n \\
\end{array}
\right),
\end{equation}
where the tildes denote that quantities updated through elimination.
The backward substition proceeds along similar lines, though it
requires the solution of linear systems with mass matrices in
dimension $d-1$.  Multiplying through each block row by
$\frac{1}{\widetilde{\nu}_{i,i}} (M^{d-1,n-i})^{-1}$ then gives, using~\eqref{eq:elevatebackward}

\begin{equation}
\label{eq:auguut}
\left(
\begin{array}{cccc|c}
I & \widetilde{\nu}^\prime_{0,1} E^{d-1,n-1,n} & \dots
& \widetilde{\nu}^\prime_{0,n} E^{d-1,0,n} & \widetilde{y}^\prime_0 \\
0 & I & \dots & \widetilde{\nu}^\prime_{1,n} E^{d-1,0,n-1} &
\widetilde{y}^\prime_1 \\
\vdots & \vdots & \ddots & \vdots \\
0 & 0 & \dots & I & \widetilde{y}^\prime_n \\
\end{array}
\right),
\end{equation}
where the primes denote quantities updated in the process.  We reflect
this in the updated $N$ matrix by scaling each row by its diagonal
entry as we proceed.
At this point, the last block of the solution is revealed,
and can be successively elevated, scaled, and subtracted from the right-hand
side to eliminate it from previous blocks.  This reveals the next-to last
block, and so-on.  We summarize this discussion in
Algorithm~\ref{alg:elim}. 

\begin{algorithm}
  \label{alg:elim}
  \caption{Block-wise Gaussian elimination for solving $M^{d,n} x = y$}
  \begin{algorithmic}
  \REQUIRE Input vector $y$ 
  \ENSURE On output, $y$ is overwritten with $(M^{d,n})^{-1} y$
  \STATE Initialize coefficient matrix $N_{a,b} := \frac{\binom{n}{a}\binom{n}{b}}{\binom{2n+d+1}{a+b}\left(2n+d\right)}$
  \FOR[Forward elimination]{$a := 0$  \TO $n$}
  \STATE $z \gets y_a$
  \FOR{$b := a+1$ \TO $n$}
  \STATE $z \gets (E^{n-1,d-b+1})^T z$
  \STATE $y_b \gets y_b - \frac{N_{b,a}}{N_{a,a}} z$
  \FOR[Elimination on $N$]{$c := a$ \TO $n$}
  \STATE $N_{b,c} \gets N_{b,c} - \frac{N_{b,a}N_{a,c}}{N_{a,a}}$
  \ENDFOR
  \ENDFOR
  \ENDFOR
  \FOR[Lower-dimensional inversion]{$a := 0$ \TO $n$}
  \STATE $y_a \gets \frac{1}{N_{a,a}} \left( M^{d-1,n-a,n-a}
  \right)^{-1} y_a$
  \FOR{$b:= a$ \TO $n$}
  \STATE $N_{b,a} \gets \frac{N_{b,a}}{N_{a,a}}$
  \ENDFOR
  \ENDFOR
  \FOR[Backward elimination]{$a:= n$ \TO $0$}
  \STATE $z \gets y_a$
  \FOR{$b:= a-1$ \TO $0$}
  \STATE $z = E^{d-1,n-b} z$
  \STATE $y_b \gets y_b - N_{b,a} z$
  \ENDFOR
  \ENDFOR
  \end{algorithmic}
\end{algorithm}

Since we will need to solve many linear systems with the same element
mass matrix, it makes sense to extend our elimination algorithm into
a reusable factorization.  We will derive a blockwise $L D L^T$
factorization of the element matrix, very much along the lines of the
standard factorizatin~\cite{strang}.

Let $N^{d,n}$ be the matrix of coefficients given in~\eqref{eq:nu}.
Suppose that we have its $LDL^T$ factorization
\begin{equation}
N^{d,n} = L_N^{d,n} D_N^{d,n} \left( L_N^{d,n} \right)^t,
\end{equation}
with
$\ell_{ij}$ and $d_{ii}$ the entries of $L_N^{d,n}$ and $D_N^{d,n}$, respectively.
We also define $U_N^{1,n} = D_N^{1,n} \left( L_N^{1,n} \right)^t$ with
$u_{ij} = d_{ii} \ell_{ji}$

Then, we can use the block matrix
\[
\widetilde{L}^{0}
=
\begin{pmatrix}
I & 0 & \dots & 0 \\
-\ell_{10} \left( E^{d-1,n-1,n} \right)^T & I & \dots & 0 \\
-\ell_{20} \left( E^{d-1,n-2,n} \right)^T & 0 & \dots & 0 \\
\vdots & \vdots & \ddots & \vdots \\
-\ell_{n0} \left( E^{d-1,0,n} \right)^T & 0 & \dots & I
\end{pmatrix}
\]
to act on $M^{d,n}$ to produce zeros below the diagonal in the first
block of columns. Similarly, we act on $\widetilde{L}^0 M^{d,n}$ with
\[
\widetilde{L}^1
= \begin{pmatrix}
I & 0 &\dots & 0 \\
0 & I & \dots & 0 \\
0 & -\ell_{21} \left( E^{d-1,n-2,n-1} \right)^T & \dots & 0 \\
\vdots & \vdots & \ddots & \vdots \\
0 & -\ell_{n1} \left( E^{d-1,0,n-1} \right)^T & \dots & I 
\end{pmatrix}
\]
to introduce zeros below the diagonal in the second block of columns.
Indeed, we have a sequence of block matrices $\widetilde{E}^{k}$ for 
$0 \leq k < n$ such that  $\widetilde{L}^k_{ij}$ is 
$P^{d-1}_{n-i} \times P^{d-1}_{n-j}$ with
\[
\widetilde{L}^k_{ij} =
\begin{cases}
I & \text{for $i = j$} \\
0 & \text{for $i \neq j$ and $j \neq k$} \\
0 & \text{for $i < j$ and $j = k$} \\
-\ell_{ij} \left( E^{d-1,n-i,n-j} \right)^T & \text{for $i > j$ and
  $j=k$} \\
\end{cases}
\]

Then, in fact, we have that
\[
\widetilde{L}^{n-1} \widetilde{L}^{n-2} \dots
\widetilde{L}^0 M^{d,n}
=
\begin{pmatrix}
u_{00} M^{d-1,n,n} & u_{01} M^{d-1,n,n-1} & \dots & u_{0n} M^{d-1,n,0}
\\
0 & u_{11} M^{d-1,n-1,n-1} & \dots & u_{1n} M^{d-1,n-1,0} \\
\vdots & \vdots & \ddots & \vdots \\
0 & 0 & \dots & u_{nn} M^{d-1,0,0}
\end{pmatrix}
\]

Much as with elementary row matrices for classic $LU$ factorization,
we can invert each of these
$\widetilde{L}^k$ matrices simply by flipping the sign of the
multiplier, so that
\[
\left( \widetilde{L}^k \right)^{-1}_{ij} =
\begin{cases}
I & \text{for $i = j$} \\
0 & \text{for $i \neq j$ and $j \neq k$} \\
0 & \text{for $i < j$ and $j = k$} \\
\ell_{ij} \left( E^{d-1,n-i,n-j} \right)^T & \text{for $i > j$ and
  $j=k$} \\
\end{cases}.
\]
Then, we define $L^{d,n}$ to be the inverse of these products
\begin{equation}
  \label{eq:Ldn}
L^{d,n} = 
\left( 
\widetilde{L}^{n-1} \widetilde{L}^{n-2} \dots \widetilde{L}^0
\right)^{-1}
= \left( \widetilde{L}^0 \right)^{-1}
\left( \widetilde{L}^1 \right)^{-1}
\dots
\left( \widetilde{L}^{n-1} \right)^{-1}
\end{equation}
so that $\left( L^{d,n} \right)^{-1} M^{d,n} \equiv U^{d,n}$ is block
upper triangular.   Like standard factorization, we can
also multiply the elimination matrices together so that
\begin{equation}
  \label{eq:Ldinv}
\left( L^{d,n} \right)^{-1}
=
\begin{pmatrix}
I & 0 &\dots & 0 \\
-\ell_{10} \left(E^{d-1,n-1,n}\right)^t & I & \dots & 0 \\
-\ell_{20} \left(E^{d-1,n-2,n}\right)^t & -\ell_{21} \left(E^{d-1,n-2,n-1}\right)^t & \dots & 0 \\
\vdots & \vdots & \ddots & \vdots \\
-\ell_{n0} \left(E^{d-1,0,n}\right)^t & -\ell_{n1} \left(E^{d-1,0,n-1}\right)^t & \dots & I 
\end{pmatrix}.
\end{equation}
Moreover, we can turn the block upper triangular matrix into a block
diagonal one times the transpose of $L^{d,n}$ giving a kind of block
$LDL^T$ factorization.  We factor out the pivot blocks from
each row, using~\eqref{eq:elimMtrick} so that
\[
U^{d,n} =
\begin{pmatrix}
d_{00} M^{d-1,n,n} & 0 & \dots & 0 \\
0 & d_{11} M^{d-1,n-1,n-1} & \dots & 0 \\
\vdots & \vdots & \ddots & \vdots \\
0 & 0 & \dots & d_{nn} M^{d-1,0,0}
\end{pmatrix}
\begin{pmatrix}
I & \ell_{10} E^{d-1,n-1,n} & \dots & \ell_{n0} E^{d-1,0,n} \\
0 & I & \dots & \ell_{n1} E^{d-1,0,n} \\
\vdots & \vdots & \ddots & \vdots \\
0 & 0 & \dots & I
\end{pmatrix}.
\]
The factor on the right is just $\left( L^{d,n} \right)^T$.

We introduce the block-diagonal matrix $\Delta^{d,n}$ by
\begin{equation}
  \label{eq:Delta}
\Delta_{ii} = d_{ii} M^{d-1,n-i}.
\end{equation}
Our discussion has established:
\begin{theorem}
\label{thm:ldlt}
The Bernstein mass matrix $M^{d,n}$ admits the block factorization
\begin{equation}
  M^{d,n} =
  L^{d,n} \Delta^{d,n} \left( L^{d,n} \right)^T.
\end{equation}
\end{theorem}

We can apply the decomposition inductively down spatial dimension, so
that each of 
the blocks in $\Delta^{d,n}$ can be also factored according to
Theorem~\ref{thm:ldlt}.  This fully expresses any mass matrix as a
diagonal matrix sandwiched in between sequences of sparse unit
triangular matrices.

So, computing the $LDL^T$ factorization of $M^{d,n}$ requires
computing the $LDL^T$ factorization of the one-dimensional coefficient
matrix $N^{d,n}$.  Supposing we use standard direct method such as
Cholesky factorization to solve
the one-dimensional mass matrices in the base case, we will have a
start-up cost of factoring $n+1$ matrices of size no larger than
$n+1$.  With Cholesky, this is a $\mathcal{O}(n^4)$ process, although
since the one-dimensional matrices factor into 
into Hankel matrices pre- and post-multiplied by diagonal matrices,
one could use Levinson's or Bareiss'
algorithm~\cite{bareiss1969numerical,levinson1947wiener} to obtain a
merely $\mathcal{O}(n^3)$ startup phase.

\begin{algorithm}
  \label{alg:ldlt}
  \caption{Mass inversion via block-recursive $LDL^T$ factorization
    for $d \geq 2$}
\begin{algorithmic}
  \REQUIRE $N^{d,n}$ factored as  $N^{d,n} = L D L^T$
  \REQUIRE Input vector $y$
  \ENSURE On output, $x = (M^{d,n})^{-1} y$
  \STATE Initialize vector $x \gets 0$
  \FOR[Apply $(L^{d,n})^{-1}$ to $y$, store in $x$]{$a := 0$  \TO $n$} 
  \STATE $z \gets y_a$
  \FOR{$b := a+1$ \TO $n$}
  \STATE $z \gets \left( E^{d-1,n-b+1} \right)^t$
  \STATE $x_b \gets x_b - L_{b,a} z$
  \ENDFOR
  \ENDFOR
  \FOR[Overwrite $x$ with $(\Delta^{d,n})^{-1} x$]{$a := 0$ \TO $n$}
      \STATE $x_a \gets \frac{1}{D_{a,a}} M^{d,n-a} x_a$
  \ENDFOR
  \FOR[Overwrite $x$ with $(L^{n,d})^{-T}x$]{$a := n$ \TO $0$}
    \STATE $z \gets x_a$
    \FOR {$b := a-1$ \TO $0$}
    \STATE $z \gets E^{d-1,n-b} z$
    \STATE $x_b \gets x_b - L_{b, a} z$
    \ENDFOR
  \ENDFOR
\end{algorithmic}
\end{algorithm}

Now, we also consider the cost of solving a linear system using the
block factorization, pseudocode for which is presented in Algorithm~\ref{alg:ldlt}.  In two dimensions, one must apply the inverse of
$L^{2,n}$, followed by the inverse of $\Delta^{2,n}$, accomplished by
triangular solves using pre-factored one-dimensional mass matrices,
and the inverse of $(L^{2,n})^T$.  In fact, the action of applying
$(L^{2,n})^{-1}$ requires exactly the same process as described above
for block Gaussian elimination, except the arithmetic on the $\nu$
values is handled in preprocessing.
That is, for each block $y_j$, we will need to compute
$\ell_{ij} (E^{1,j-i,j})^T y_j$ for $1 \leq i \leq n-j-1$ and
accumulate scalings of these vectors into corresponding blocks of the
result.  Since these elevations are needed for each $i$, it is
helpful to reuse these results.  Applying
$(L^{2,n})^{-1}$ then requires applying
$E^{1,i-j}$ for all valid $i$ and $j$, together with all of the
axpy operations.  Since the one-dimensional elevation into degree $i$
has $2(i+1)$ nonzeros in it, the required elevations required cost
\begin{equation}
  \label{eq:elev2count}
\sum_{i=1}^n \sum_{j=1}^{i-1} 2(j+1) = \frac{n(n^2+3n-4)}{3},
\end{equation}
operations, which is $\mathcal{O}(n^3)$, 
and we also have a comparable number of operations for the axpy-like
operations to accumulate the result.  A similar discussion shows that 
applying $(L^{2,n})^{-T}$ requires the same number of operations.
Between these stages, one must invert the lower-dimensional mass
matrices using the pre-computed Cholesky factorizations and perform
the scalings to apply $\Delta^{-1}$.  Since a 
pair of $m \times m$ triangular solves costs $m(m+1)$
operations, the total cost of the one-dimensional mass inversions is
\[
\sum_{i=0}^{n} (i+1)(i+2) = \frac{(n+1)(n+2)(n+3)}{3},
\]
together with the lower-order term for scalings
\[
\sum_{i=0}^{n} P^1_i = \sum_{i=0}^{n} (i+1) = \frac{(n+1)(n+2)}{2}.
\]
So, the whole three-stage process is $\mathcal{O}(n^3)$ per element.

In dimension $d > 2$, we may proceed inductively in space dimension to
show that Algorithm~\ref{alg:ldlt} requires, after start-up,
$\mathcal{O}(n^{d+1})$ operations.  The application of $\Delta^{-1}$
will always require $n+1$ inversions of $(d-1)$-dimensional mass
matrices , each of which costs $\mathcal{O}(n^{d})$ operations by the
induction hypothesis.  Inverting $\Delta^{d,n}$ onto a vector will cost
$\mathcal{O}(n^{d+1})$ operations for all $n$ and $d$.  To see that a
similar complexity holds for applying the inverses of $L^{d,n}$ and
its transpose, one can simply replace the summand
in~\eqref{eq:elev2count} with $2 P^{d-1}_j$ and execute the sum.  To
conclude,
\begin{theorem}
  Algorithm~\ref{alg:ldlt} applies the inverse of $M^{d,n}$ to an
  arbitrary vector in $\mathcal{O}(n^{d+1})$ operations.
\end{theorem}

\section{Numerical results}

\subsection{Mass inversion}
\begin{figure}
  \begin{center}
    \includegraphics[width=4.0in]{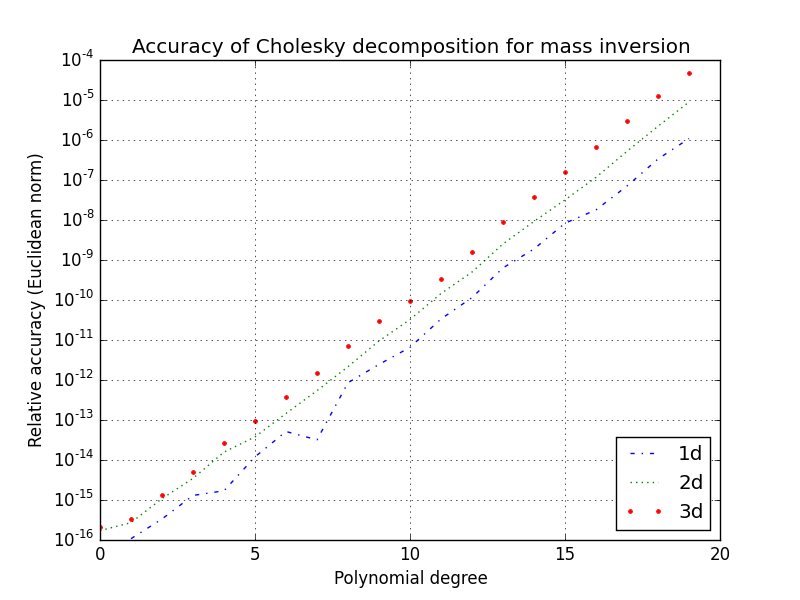}
    \caption{Relative accuracy of solving linear systems with mass
      matrices of various degrees using Cholesky decomposition.}  
  \label{fig:cholacc}
  \end{center}
\end{figure}

Because of Corollary~\ref{cor:badguy}, we must pay special attention
to the accuracy with which linear systems involving the mass matrix
are computed.  We began with Cholesky decomposition as a baseline.
For degrees one through twenty in one, two, and three space dimension,
we explicitly formed the reference mass matrix in Python and used the
\texttt{scipy}~\citep{jones2001scipy} interface LAPACK to form
the Cholesky decomposition.  Then, we chose several
random vectors to be sample solutions and formed the right-hand side by 
direct matrix-vector multiplication.  In Figure~\ref{fig:cholacc}, we
plot the relative accuracy of a function of degree in each space
dimension.  Although we observe expontial growth in the error (fully
expected in light of Corollary~\ref{cor:badguy}), we see that we still
obtain at least ten digits of relative accuracy up to degree ten.

\begin{figure}
  \begin{center}
  \subfigure[Accuracy obtained by iterating until a 
    residual tolerance of $10^{-12}$.]{%
    \includegraphics[width=0.4\textwidth]{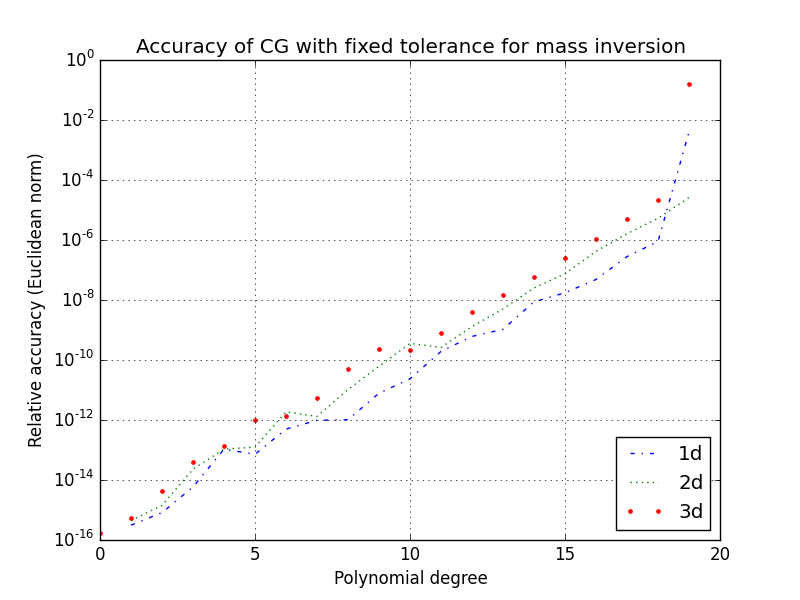}\label{acc}
  }
  \hspace{5mm}
  \subfigure[CG iterations required to solve $M^{d,n} x = y$ to a tolerance of $10^{-12}$.
  ]{%
    \includegraphics[width=0.4\textwidth]{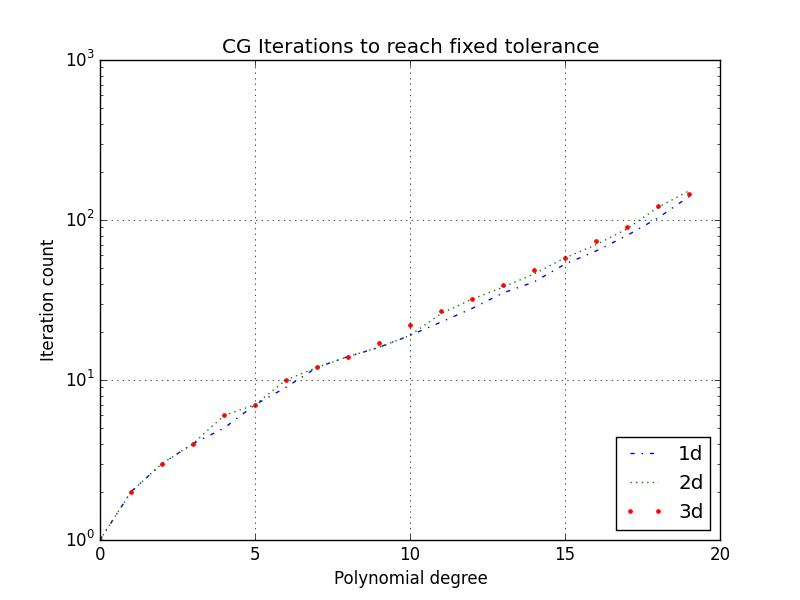}\label{its}
  }
  \caption{%
Accuracy obtained solving mass matrix system using conjugate gradient
iteration in one, two, and three space dimensions.  
  }
      \label{fig:cgtol}
  \end{center}
\end{figure}

\begin{figure}
  \label{fig:cgits}
  \begin{center}
    \includegraphics[width=4.0in]{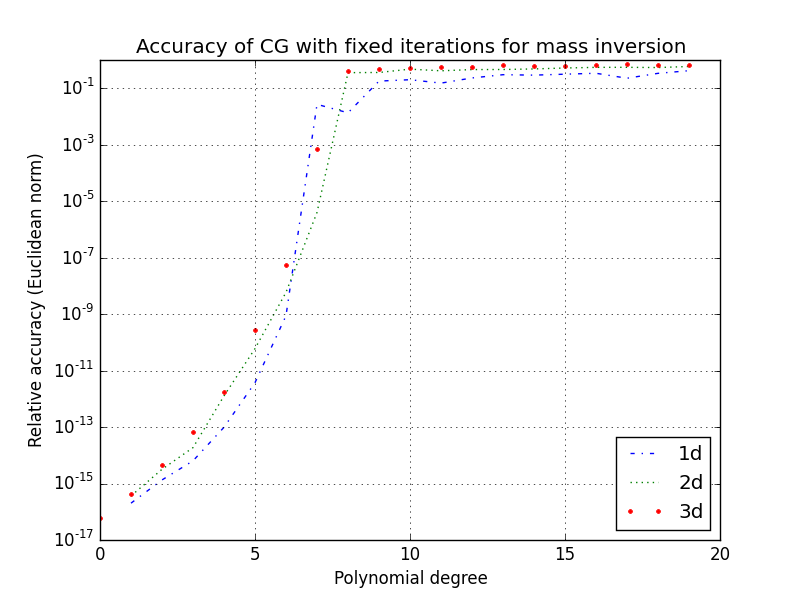}
    \caption{Relative accuracy of solving $M^{d,n}x = y$ using
      exactly $n+1$ CG iterations.}
  \end{center}
\end{figure}

Second, we also attempt to solve the linear system using conjugate
gradients.  We again used systems  with random solution,
and both letting CG run to a relative residual tolerance of
$10^{-12}$ and also stopping after $n+1$ iterations in light of
Corollary~\ref{cor:wonthappen}.  We display the results of a fixed
tolerance in Figure~\ref{fig:cgtol}.  Figure~\ref{acc},
shows the actual accuracy obtained for each polynomial degree and
Figure~\ref{its} gives the actual iteration count required.  Like
Cholesky factorization, this approach gives nearly ten-digit accuracy
up to degree ten polynomials.  On the other hand,
Figure~\ref{fig:cgits} shows that accuracy degrades markedly when
only $n+1$ iterations are used.

\begin{figure}
  \begin{center}
    \includegraphics[width=4.0in]{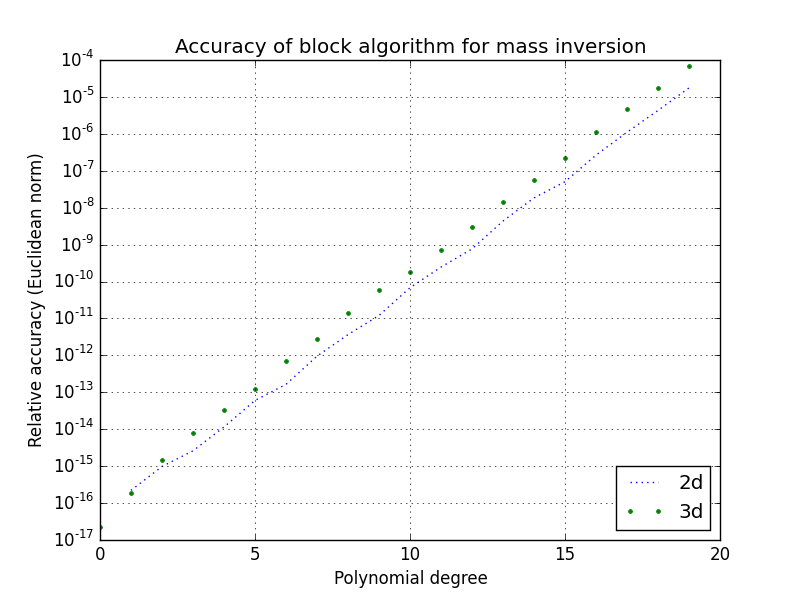}
    \caption{Relative accuracy of solving linear systems with mass
      matrices of various degrees using one level of the block
      algorithm with Cholesky factorization for lower-dimensional
      matrices.}
      \label{fig:blockacc}
  \end{center}
\end{figure}

Finally, our block algorithm gives accuracy comparable to that of
Cholesky factorization.  Our two-dimensional implementation of
Algorithm~\ref{alg:ldlt} uses Cholesky factorizations of the
one-dimensional mass matrices.  Rather than full recursion, our
three-dimensional implementation uses Cholesky factorization of the
two-dimensional matrices.  At any rate, Figure~\ref{fig:blockacc}
shows, when compared to Figure~\ref{fig:cholacc}, that we lose very
little additional accuracy over Cholesky factorization.  Whether
replacing the one-dimensional solver with a specialized method
for totally positive matrices~\cite{koev2007accurate} would also give
high accuracy for the higher-dimensional problems will be the subject
of future investigation.

\subsection{Timing for first-order acoustics}
\begin{figure}
  \label{fig:2dtimings}
  \begin{center}
    \includegraphics[width=4.0in]{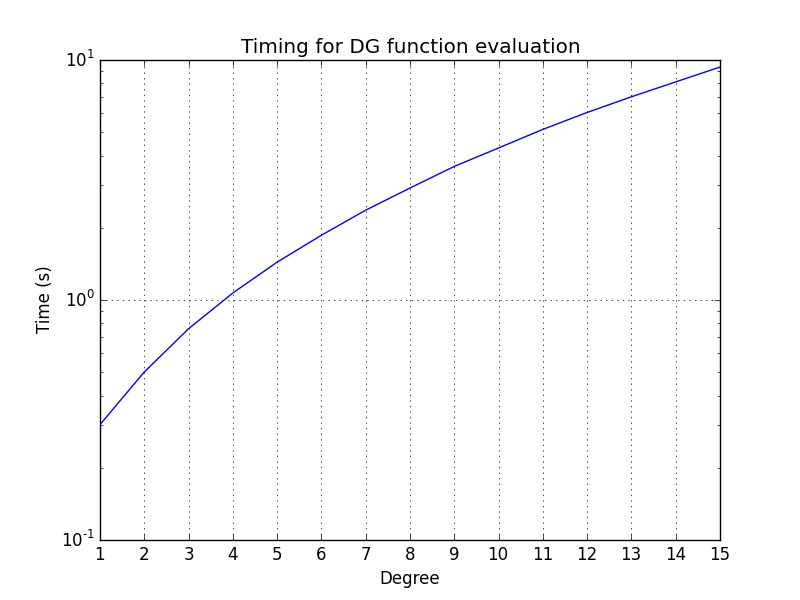}
    \caption{Timing of a DG function evaluation for various polynomial
      degrees on a $32 \times 32$ mesh.}
  \end{center}
\end{figure}

We fixed a $32 \times 32$ square mesh subdivided into right triangles
and computed the time to perform the DG function evaluation (including
mass matrix inversion) at various polynomial degrees.  We used the
mesh from DOLFIN~\citep{dolfin} and wrote the Bernstein polynomial algorithms
in Cython~\citep{behnelcython}.
With an $\mathcal{O}(n^3)$ complexity for two-dimensional problems, we
expect a doubling of the polynomial degree to produce an eightfold
increase in run-time.  In Figure~\ref{fig:2dtimings}, though, we see
even better results.  In fact, a least-squares fit of the log-log data
in this table from degrees five to fifteen gives a very near fit with
a slope of less than two (about 1.7) rather than three.  Since small
calculations tend to run at lower flop rates, it is possible that we
are far from the asymptotical regime predicted by our operation counts.

\section{Conclusions and Future Work}
Bernstein polynomials admit optimal-complexity algorithms for
discontinuous Galerkin methods for conservation laws.  The dense
element mass matrices might, at first blush, seem to prevent this,
but their dimensionally recursive block structure and other interesting properties, lead to an efficient blockwise factorzation.  Despite the large
condition numbers, our current algorithms seem sufficient to deliver
reasonable accuracy at moderate polynomial orders.

On the other hand, these results still leave much room for future
investigation.  First, it makes sense to explore the possibilities of
slope limiting in the Bernstein basis.  Second, while our 
mass inversion algorithm is sufficient for moderate order, it may be
possible to construct a different algorithm that maintains the low
complexity while giving higher relative accuracy, enabling very high
approximation orders.  Perhaps such algorithms will either utilize the techniques in~\citep{koev2007accurate} internally, or else extend them somehow.
Finally, our new algorithm, while of optimal compexity, is quite
intricate to implement and still is not well-tuned for high
performance.  Finding ways to make these algorithms more performant
will have important practical benefits.

\bibliographystyle{siam} 
\bibliography{bibliography}
\end{document}